\documentclass[12pt]{amsart}

\title[Infinite chess position with value $\omega^4$]{A position in infinite chess with game value $\omega^4$}

\author[Evans]{C.~D.~A.~Evans}
 \address[C.~D.~A.~Evans]
         {Program in Philosophy,
          The Graduate Center of The City University of New York,
          365 Fifth Avenue, New York, NY 10016}
 \email{c.alexander.evans@gmail.com}
 \urladdr{http://please.provide}
 \thanks{The first author holds the chess title of National Master (USA)}

\author[Hamkins]{Joel David Hamkins}
 \address[J.~D.~Hamkins]
         {Philosophy, New York University, 5 Washington Place New York, New York 10003 \&
         Mathematics, The Graduate Center of The City University of New York,
         365 Fifth Avenue, New York, NY 10016 \&
         Mathematics, College of Staten Island of CUNY, Staten Island, NY 10314}
 \email{jhamkins@gc.cuny.edu}
 \urladdr{http://jdh.hamkins.org}
 \thanks{The research of the second author has been supported in part by NSF grant
  DMS-0800762, PSC-CUNY grant 64732-00-42 and Simons Foundation grant 209252. In addition, part of this research was undertaken whilst the second author was a visiting fellow at the Isaac Newton Institute for Mathematical Sciences in the programme `Mathematical, Foundational and Computational Aspects of the Higher Infinite' (HIF). Commentary concerning this paper can be made at \href{http://jdh.hamkins.org/a-position-in-infinite-chess-with-game-value-omega-to-the-4}{jdh.hamkins.org/a-position-in-infinite-chess-with-game-value-omega-to-the-4}.}

\author[Perlmutter]{Norman Lewis Perlmutter}
 \address[N.~Perlmutter]
         {Mathematics, LaGuardia Community College of CUNY}
 \email{drplaguardia@gmail.com}
 \urladdr{http://boolesrings.org/perlmutter}

\newtheorem{theorem}{Theorem}

\newtheorem*{questions*}{Questions}
\newtheorem*{mainquestion*}{Main Question} 
\newtheorem*{openquestion*}{Open Question} 

\newcommand{\QED}{\end{proof}}

\def\proclaim[#1]{{\bf #1}}
\def\BF#1.{{\bf #1.}}

\newcommand{\Z}{{\mathbb Z}}

\newcommand{\omegaoneCh}{\omega_1^\Ch}
\newcommand{\omegaoneChi}{\omega_1^{\baselineskip=0pt\vtop to 7pt{\hbox{$\scriptstyle\Ch$}\vskip-1pt\hbox{\,$\scriptscriptstyle\sim$}}}}
\newcommand{\omegaoneChc}{\omega_1^{\Ch,c}}

\newcommand{\omegaoneChthreei}{\omega_1^{\baselineskip=0pt\vtop to 7pt{\hbox{$\scriptstyle\Ch_3$}\vskip-1.5pt\hbox{\,$\scriptscriptstyle\sim$}}}}
\newcommand{\omegaoneChthreec}{\omega_1^{\Ch_3,c}}
\newcommand{\smalllt}{\mathrel{\mathchoice{\raise2pt\hbox{$\scriptstyle<$}}{\raise1pt\hbox{$\scriptstyle<$}}{\raise0pt\hbox{$\scriptscriptstyle<$}}{\scriptscriptstyle<}}}
\newcommand{\smallleq}{\mathrel{\mathchoice{\raise2pt\hbox{$\scriptstyle\leq$}}{\raise1pt\hbox{$\scriptstyle\leq$}}{\raise1pt\hbox{$\scriptscriptstyle\leq$}}{\scriptscriptstyle\leq}}}

\newcommand{\boolval}[1]{\mathopen{\lbrack\!\lbrack}\,#1\,\mathclose{\rbrack\!\rbrack}}
\def\[#1]{\boolval{#1}}
\newbox\gnBoxA
\newdimen\gnCornerHgt
\setbox\gnBoxA=\hbox{\tiny$\ulcorner$}
\global\gnCornerHgt=\ht\gnBoxA
\newdimen\gnArgHgt
\def\gcode #1{%
\setbox\gnBoxA=\hbox{$#1$}%
\gnArgHgt=\ht\gnBoxA%
\ifnum     \gnArgHgt<\gnCornerHgt \gnArgHgt=0pt%
\else \advance \gnArgHgt by -\gnCornerHgt%
\fi \raise\gnArgHgt\hbox{\tiny$\ulcorner$} \box\gnBoxA %
\raise\gnArgHgt\hbox{\tiny$\urcorner$}}
\newcommand{\UnderTilde}[1]{{\setbox1=\hbox{$#1$}\baselineskip=0pt\vtop{\hbox{$#1$}\hbox to\wd1{\hfil$\sim$\hfil}}}{}}
\newcommand{\Undertilde}[1]{{\setbox1=\hbox{$#1$}\baselineskip=0pt\vtop{\hbox{$#1$}\hbox to\wd1{\hfil$\scriptstyle\sim$\hfil}}}{}}
\newcommand{\undertilde}[1]{{\setbox1=\hbox{$#1$}\baselineskip=0pt\vtop{\hbox{$#1$}\hbox to\wd1{\hfil$\scriptscriptstyle\sim$\hfil}}}{}}
\newcommand{\UnderdTilde}[1]{{\setbox1=\hbox{$#1$}\baselineskip=0pt\vtop{\hbox{$#1$}\hbox to\wd1{\hfil$\approx$\hfil}}}{}}
\newcommand{\Underdtilde}[1]{{\setbox1=\hbox{$#1$}\baselineskip=0pt\vtop{\hbox{$#1$}\hbox to\wd1{\hfil\scriptsize$\approx$\hfil}}}{}}

\def\<#1>{\left\langle#1\right\rangle}

\newcommand{\Ch}{{\mathfrak{Ch}}}

\newcommand{\omegaCK}{{\omega_1^{\hbox{\tiny\sc CK}}}}
\newcommand{\cell}[1]{\boxit{\hbox to 17pt{\strut\hfil$#1$\hfil}}}
\newcommand{\head}[2]{\lower2pt\vbox{\hbox{\strut\footnotesize\it\hskip3pt#2}\boxit{\cell#1}}}
\newcommand{\boxit}[1]{\setbox4=\hbox{\kern2pt#1\kern2pt}\hbox{\vrule\vbox{\hrule\kern2pt\box4\kern2pt\hrule}\vrule}}
\newcommand{\Col}[3]{\hbox{\vbox{\baselineskip=0pt\parskip=0pt\cell#1\cell#2\cell#3}}}
\newcommand{\tapenames}{\raise 5pt\vbox to .7in{\hbox to .8in{\it\hfill input: \strut}\vfill\hbox to
.8in{\it\hfill scratch: \strut}\vfill\hbox to .8in{\it\hfill output: \strut}}}
\newcommand{\Head}[4]{\lower2pt\vbox{\hbox to25pt{\strut\footnotesize\it\hfill#4\hfill}\boxit{\Col#1#2#3}}}
\newcommand{\Dots}{\raise 5pt\vbox to .7in{\hbox{\ $\cdots$\strut}\vfill\hbox{\ $\cdots$\strut}\vfill\hbox{\
$\cdots$\strut}}}
\newcommand{\df}{\it} 
\usepackage{chessboard}
\newcommand\forcedmove{{\tiny\bf$\Box$}}
\usepackage[hidelinks]{hyperref}
\usepackage{latexsym,amsfonts,amsmath,amssymb}
\usepackage{tikz} 
\usetikzlibrary{arrows}
\usepackage{diagrams}\diagramstyle[tight,centredisplay,textflow]

\begin{document}

\begin{abstract}
  We present a position in infinite chess exhibiting an ordinal game value of $\omega^4$, thereby
  improving on the previously largest-known values of $\omega^3$ and $\omega^3\cdot 4$.
\end{abstract}

\maketitle

Infinite chess is chess played on an infinite chessboard, a field of alternating black and white squares arranged like the integer lattice $\Z\times\Z$, upon which the familiar chess characters---kings, queens, rooks, bishops, knights and pawns---move about with the aim to capture the opposing king. Since checkmate, when it occurs, does so at a finite stage of play, it follows that infinite chess is what is known technically as an open game, and it is therefore subject to the theory of transfinite ordinal game values, which is applicable in any open game. In particular, a designated player in a game of infinite chess has a winning strategy from a given position if and only if that position has an ordinal game value, and it is a winning strategy for that player simply to play so as to reduce value. In such a case, in an open game proceeding from a position with high ordinal game value for white, say, then although white will ultimately win the game, nevertheless the length of play is controlled by black in a way that is tightly connected with the size of the ordinal: playing an open game with value $\alpha$ for white is essentially analogous to having black count down in the ordinals from $\alpha$, with white winning when $0$ is reached, as inevitably it will be. For this reason, these games can often have a somber absurd character of play, where white slowly makes progress on a main-line attack while black, although doomed, mounts aggressive but increasingly desperate countermeasures including long bouts of prolonged forced play that distract white's progress in the main attack and stave off the inevitable checkmate a bit longer, but not indefinitely.

It is an open question how large the game values arising for positions in infinite chess can be. The largest values previously known were provided by Evans and Hamkins~\cite{EvansHamkins2014:TransfiniteGameValuesInInfiniteChess}  (the first two authors of this article), who presented a position with value $\omega^3$ and who indicated how similar positions could be expected to achieve values as high as $\omega^3\cdot 4$. They defined the {\df omega one of chess} to be the supremum of the game values that arise for positions in infinite chess. For the boldface version $\omegaoneChi$, one takes this supremum over all positions, including those with infinitely many pieces, while for the lightface version $\omegaoneCh$, one considers only the positions having only finitely many pieces; there is also a computable version $\omegaoneChc$, where one considers only computable positions, and other alternative versions of the omega one of chess correspond to each complexity class to which the position might belong. Some natural upper bounds for these values are provided by~\cite[thm.~4]{EvansHamkins2014:TransfiniteGameValuesInInfiniteChess} as follows, where $\omegaCK$ is the Church-Kleene ordinal.
$$\omegaoneCh\leq\omegaoneChc\leq\omegaCK$$
$$\omegaoneChi\leq\omega_1.$$
Evans and Hamkins~\cite[conj.~5]{EvansHamkins2014:TransfiniteGameValuesInInfiniteChess} have conjectured that these upper bounds are sharp, even though the best lower bounds we have currently are significantly smaller.\footnote{In the case of infinite three-dimensional chess, however, the corresponding upper bounds are realized. Specifically, Evans and Hamkins~\cite[thms.~10,11]{EvansHamkins2014:TransfiniteGameValuesInInfiniteChess} proved that every countable ordinal arises as the value of a game position in infinite three-dimensional chess, so that $\omegaoneChthreei=\omega_1$, and every computable ordinal arises as the value of a computable position in infinite three-dimensional chess, so that $\omegaoneChthreec=\omegaCK$.}

In this article, we shall improve the lower bound by providing a computable position with value $\omega^4$, a new record.

\begin{theorem}
 There is a position in infinite chess with game value $\omega^4$. The position is computable and involves infinitely many pieces.
\end{theorem}

\begin{figure}[h]
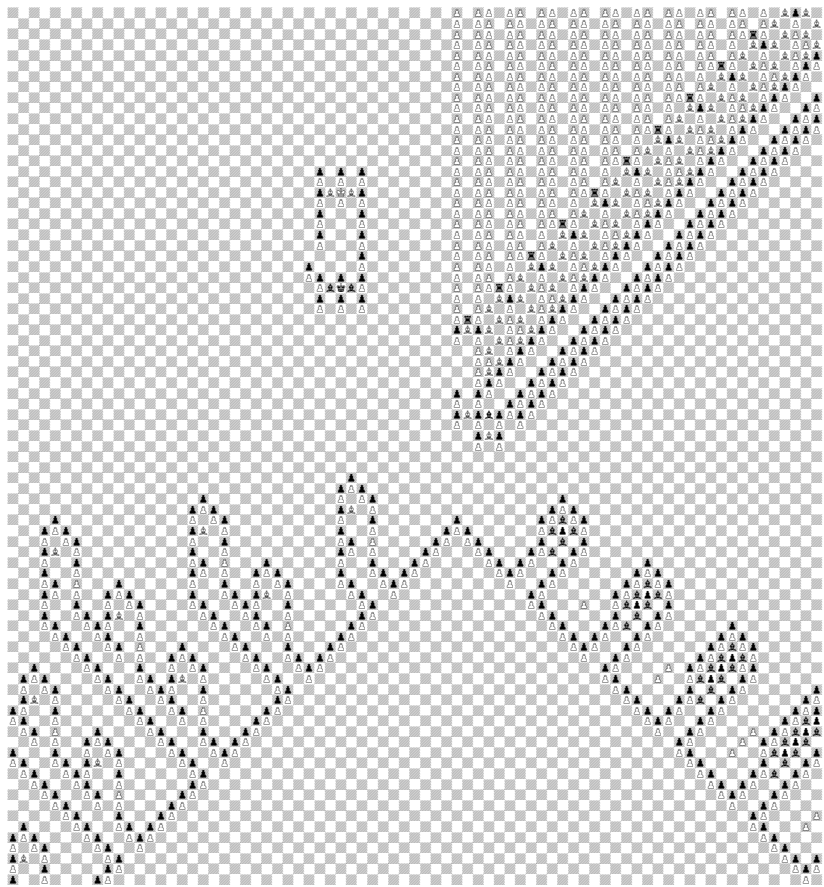

\catcode`\.=\active\def.{{1}}
\chessboard[maxfield=z43,
            boardfontsize=4pt,
            label=false,
            showmover=false,
            borderleft=false,
            borderright=false,
            bordertop=false,
            borderbottom=false,
            margin=false,
            setfen=%
../%
../%
../%
../%
../%
../%
../%
../%
../%
../%
../%
../%
../%
../%
../%
../%
../%
../%
../%
../%
../%
../%
../%
../%
../%
../%
../%
../%
../%
../%
../%
../%
../%
../%
../%
../%
../%
../%
../%
../%
../%
../%
../%
../%
../%
../%
../%
]%
\chessboard[maxfield=p43,
            boardfontsize=4pt,
            label=false,
            showmover=false,
            borderleft=false,
            borderright=false,
            bordertop=false,
            borderbottom=false,
            margin=false,
            setfen=%
../%
.../%
../%
.../%
../%
.../%
../%
.../%
../%
.../%
../%
...../%
....../%
...../%
...../%
...p.p.p../%
...P.P.P../%
...pBKBp../%
...P.P.P../%
...p...p../%
...P...P../%
...p...p../%
...P...P../%
.......p../%
..p....P../%
..Pp.p.p../%
...PbkbP../%
...p.p.p../%
...P.P.P../%
../%
../%
../%
../%
../%
../%
]%
\chessboard[maxfield=p43,
            boardfontsize=4pt,
            labelbottom=false,
            labelleft=false,
            labelfontsize=6pt,
            labelleftwidth=2ex,
            showmover=false,
            borderleft=false,
            borderright=false,
            bordertop=false,
            borderbottom=false,
            coloremphstyle=\color{green!25!black},
            coloremphstyle=\color{blue!50!black},
            coloremphstyle=\color{red!50!black},
            pgfstyle=straightmove,linewidth=1pt,
            margin=false,
            setfen=%
P.PP.PP.PP.PP.PP/%
P.PP.PP.PP.PP.PP/%
P.PP.PP.PP.PP.PP/%
P.PP.PP.PP.PP.PP/%
P.PP.PP.PP.PP.PP/%
P.PP.PP.PP.PP.PP/%
P.PP.PP.PP.PP.PP/%
P.PP.PP.PP.PP.PP/%
P.PP.PP.PP.PP.PP/%
P.PP.PP.PP.PP.PP/%
P.PP.PP.PP.PP.PP/%
P.PP.PP.PP.PP.PP/%
P.PP.PP.PP.PP.PP/%
P.PP.PP.PP.PP.PP/%
P.PP.PP.PP.PP.PP/%
P.PP.PP.PP.PP.P./%
P.PP.PP.PP.PP.PB/%
P.PP.PP.PP.PPrP./%
P.PP.PP.PP.P.BpB/%
P.PP.PP.PP.PB.P./%
P.PP.PP.PPrP.BPB/%
P.PP.PP.P.BpB.PP/%
P.PP.PP.PB.P.BPB/%
P.PP.PPrP.BPB.Pp/%
P.PP.P.BpB.PPBpP/%
P.PP.PB.P.BPBpP./%
P.PPrP.BPB.PpP../%
P.P.BpB.PPBpP..p/%
P.PB.P.BPBpP..pP/%
PrP.BPB.PpP..pPp/%
pBpB.PPBpP..pPpP/%
P.P.BPBpP..pPpP/%
..PB.PpP..pPpP/%
..PPBpP..pPpP/%
..PBpP..pPpP/%
..PpP..pPpP............/%
p.pP..pPpP............../%
P.P..pPpP............/%
pBpbpPpP.........../%
P.P.P.P............/%
..pBp................/%
..P.P.........../%
........../%
........../%
]%
\chessboard[maxfield=s43,
            boardfontsize=4pt,
            label=false,
            showmover=false,
            mover=b,
            borderleft=false,
            borderright=false,
            bordertop=false,
            borderbottom=false,
            margin=false,
            coloremphstyle=\color{blue!50!black},
            setfen=%
.PP.PP.PP.PP.P.BpB.PPBpP..pP/%
.PP.PP.PP.PP.PB.P.BPBpP..pPp/%
.PP.PP.PP.PPrP.BPB.PpP..pPpP/%
.PP.PP.PP.P.BpB.PPBpP..pPpP/%
.PP.PP.PP.PB.P.BPBpP..pPpP/%
.PP.PP.PPrP.BPB.PpP..pPpP/%
.PP.PP.P.BpB.PPBpP..pPpP/%
.PP.PP.PB.P.BPBpP..pPpP/%
.PP.PPrP.BPB.PpP..pPpP/%
.PP.P.BpB.PPBpP..pPpP/%
.PP.PB.P.BPBpP..pPpP/%
.PPrP.BPB.PpP..pPpP/%
.P.BpB.PPBpP..pPpP/%
.PB.P.BPBpP..pPpP/%
rP.BPB.PpP..pPpP/%
BpB.PPBpP..pPpP/%
.P.BPBpP..pPpP/%
BPB.PpP..pPpP/%
.PPBpP..pPpP/%
BPBpP..pPpP/%
.PpP..pPpP/%
BpP..pPpP/%
pP..pPpP/%
P..pPpP/%
..pPpP/%
.pPpP/%
pPpP........./%
PpP/%
pP/%
P/%
]
\chessboard[maxfield=v40,
            boardfontsize=4pt,
            label=false,
            showmover=false,
            borderleft=false,
            borderright=false,
            bordertop=false,
            borderbottom=false,
            margin=false,
            setfen=%
....................../%
....................../%
....................../%
..................p.../%
.................pPp../%
....p............P.Pp.../%
...pPp...........pB.P.../%
...P.Pp..........P..p.../%
...pB.P..........p..P.../%
...P..p..........Pp.P./%
...p..P..........pP.P./%
...Pp.P...p......P..p../%
...pP.P..pPp.....p..Pp./%
...P..p..P.Pp....Pp..P../%
...p..Pp.pB.P.....Pp.../%
...Pp..PpP..p......Pp.../%
....Pp..Pp..P.......Pp./%
.....Pp..Pp.P...p....P/%
......Pp..P.P..pPp..../%
..p....Pp...p..P.Pp......../%
.pPp....Pp..Pp.pB.P....../%
.P.Pp....Pp..PpP..p......../%
.pB.P.....Pp..Pp..P......./%
pP..p......Pp..Pp.P...../%
Pp..P.......Pp..P.P...../%
.Pp.P...p....Pp...p..../%
..P.P..pPp....Pp..Pp.p.../%
p...p..P.Pp....Pp..PpP..../%
Pp..Pp.pB.P.....Pp..P...../%
.Pp..PpP..p......Pp..../%
..Pp..Pp..P......pP..../%
...Pp..Pp.P.....pP..../%
....Pp..P.P....pP...../%
.....Pp...p...pP...../%
.p....Pp..Pp.pP...../%
pPp....Pp..PpP...../%
P.Pp....Pp..P...../%
pB.P.....Pp....../%
P..p.....pP...../%
p..P....pP...../%
Pp.P...pP...../%
.P.P..pP...../%
...p.pP.....%
p..PpP..../%
Pp..P..../%
.Pp....../%
.pP..../%
pP./%
]%
\chessboard[maxfield=t40,
            boardfontsize=4pt,
            label=false,
            showmover=false,
            borderleft=false,
            borderright=false,
            bordertop=false,
            borderbottom=false,
            margin=false,
            setfen=%
..................../%
..........p......../%
.........pPp......../%
.........P.Pp....../%
.........pB.P......../%
.........P..p......../%
.........p..P......p../%
.........Pp.P.....pP./%
.........pP.P....pP...../%
..p......P..p...pP....../%
.pPp.....p..Pp.pP../%
.P.Pp....Pp..PpP../%
.pB.P.....Pp..P.../%
pP..p......Pp..../%
Pp..P......pP...../%
.Pp.P.....pP...../%
..P.P....pP..../%
p...p...pP...../%
Pp..Pp.pP...../%
.Pp..PpP..../%
..Pp..P...../%
...Pp....../%
...pP..../%
..pP...../%
.pP...../%
pP..../%
P/%
./%
]%
\chessboard[maxfield=z40,
            boardfontsize=4pt,
            label=false,
            showmover=false,
            borderleft=false,
            borderright=false,
            bordertop=false,
            borderbottom=false,
            margin=false,
            setfen=%
................................../%
............................../%
............................../%
..........p.................../%
.........pPp................./%
p.......pPbPp.............../%
Pp......PbpbP.............../%
.Pp.....p.b.p............../%
..Pp...pPb.pP............./%
...Pp.pP..pP......p......./%
....PpP..pP......pPp.........../%
.....P..pP......pPbPp...../%
.......pP......pPbpbP....../%
.......Pp...P..Pbpb.p.........../%
........Pp.....p.b.pP......../%
.........Pp...pPb.pP........../%
..........Pp.pP..pP......p/%
...........PpP..pP......pP/%
............P..pP......pPb/%
..............pP....P.pPbp/%
..............Pp...P..Pbpb/%
...............Pp.....p.b./%
................Pp...pPb.p..../%
.................Pp.pP..pP/%
..................PpP..pP./%
...................P..pP../%
.....................pP.../%
.....................Pp.../%
......................Pp../%
.......................Pp./%
........................Pp/%
.........................P/%
]%
\chessboard[maxfield=i40,
            boardfontsize=4pt,
            label=false,
            showmover=false,
            borderleft=false,
            borderright=false,
            bordertop=false,
            borderbottom=false,
            margin=false,
            setfen=%
./%
./%
./%
..../%
../%
./%
/%
/%
/%
/%
/%
/%
/%
../%
......../%
p...../%
Pp....../%
bPp../%
pbP../%
bPp...../%
.pP/%
pP......p/%
P......pPp/%
......pPpPp./%
.....pPbpbP./%
..P.pPbpb.p/%
.P.pPbpb.pP/%
P..Pbpb.pP/%
...p.b.pP./%
..pPb.pP./%
.pP..pP../%
pP..pP....P.p/%
P..pP....P.pP/%
..pP....P.pPb/%
..Pp...P..Pbp/%
...Pp.....p.b/%
....Pp...pPb./%
.....Pp.pP..p/%
......PpP..pP/%
.......P..pP/%
]%
\caption{A position in infinite chess with value $\omega^4$}\label{Figure:MainPositionValueOmega^4}
\end{figure}

The position, shown in figure~\ref{Figure:MainPositionValueOmega^4}, is analyzed in section~\ref{Section.APositionWithValueOmega^4}, where we argue that the game value is $\omega^4$.

This article continues the work of the first two authors in~\cite{EvansHamkins2014:TransfiniteGameValuesInInfiniteChess}, inspired by the MathOverflow questions of Johan W\"astlund~\cite{MO63423Wastlund:CheckmateInOmegaMoves} and Richard Stanley~\cite{MO27967Stanley:DecidabilityOfInfiniteChess}, which had also led to the second author's joint work~\cite{BrumleveHamkinsSchlicht2012:TheMateInNProblemOfInfiniteChessIsDecidable}.

\section{Transfinite game values in infinite chess}\label{Section.TransfiniteGameValues}

The reader is encouraged to consult~\cite[\S1]{EvansHamkins2014:TransfiniteGameValuesInInfiniteChess} for
a general discussion of game values in the context of infinite chess. We briefly recount here some of the
main points. A two-player infinite game of perfect information---one in which the players alternately
make their moves in turn, with the winner determined by the sequence $m_0,m_1,\ldots$ of these moves---is
{\df open} for a given player if whenever a play of the game is won for that player, then it is won after only
finitely many moves. This is equivalent to saying that the collection of all sequences of play that are
winning for that player form an open set in the product topology in the space of all possible plays of
the game.

In any open game, one may define the ordinal game valuation for a designated player on the set of possible positions in the game as follows. Since we are thinking of chess, let us call the players white and black and define the game value from the perspective of white (a dual valuation applies from the perspective of black).
\begin{itemize}
 \item If the position shows that white has already won, then the position has value $0$.
 \item If it is white's turn to play from position $p$, then the value of $p$ is $\alpha+1$, if $\alpha$ is minimal such that white may make a move from $p$ to a new position having value $\alpha$.
 \item If it is black's turn to play from a position $p$, then the value of $p$ is defined only if every legal move by black from $p$ results in a new position with a value, and in this case the value of $p$ is the supremum of the values of those resulting positions.
\end{itemize}
The game valuation is the minimal assignment of ordinal values that obeys these recursive rules; some positions may be left without any value assigned.

The {\bf fundamental observation of game values} is that if a position has an ordinal value, then white
has a winning strategy for play proceeding from that position; and if it lacks a value, then black has a
strategy to prevent white from winning. This observation amounts to a proof of the Gale-Stewart theorem
\cite{GaleStewart1953:InfiniteGamesWithPerfectInformation}, asserting that all open games are determined,
and it is surprisingly easy to prove. Namely, if a position has a value, then since black cannot play so
as to increase value and white can play so as reduce it (unless it is already zero), it follows that the value-reducing strategy will lead to a strictly descending sequence of ordinals during the course of play. Since the ordinals are
well-ordered, there is no infinite descending sequence of ordinals, and so the value-reducing strategy
must eventually arrive at the value $0$, which is a win for white. Therefore, white has a winning
strategy in the game proceeding from any position with a value. Conversely, if a position does not have a
value, then white cannot play so as to arrive at a position with a value, and black can play so as to
maintain the fact that it has no value, which means that black has a strategy, the maintaining strategy,
to prevent white from winning; in other words, black can force either a draw or a win for black. Black
will have a winning strategy just in case the position has an ordinal value when defined with black as
the designated player. So in an open game like chess that allows for drawn results, the situation is that
if a position has a game value for either player, then that player has a winning strategy from that
position, and otherwise, both players have strategies to force a draw. Note that in infinite chess, we abandon the fifty-move tournament rule of ordinary $8\times 8$ chess as an unnatural restriction, and the threefold-repetition rule is unnecessary, because a game is regarded as drawn if it lasts forever without checkmate; draws by stalemate, however, can still occur in infinite chess.

Let us illustrate the concrete meaning of ordinal game values in the case of various small ordinals. It is easy to see, for example, that a position has finite value $n$ just in case white can force a win in $n$ moves, but not quicker than this (counting the moves only of white). A position has value at least $\omega$ just in case it is black to play, and white has a winning strategy that will force a win in finitely many moves, but black may announce with his first move a lower bound for the length of play. A position has value at least $\omega\cdot n$, if white has a winning strategy, but black can make $n$ such announcements, where each announcement is the number of moves until the next announcement. A position has value at least $\omega^2$, if white has a winning strategy, but at the start of play black can make a large announcement, a number $n$, which the number of announcements of the previous type that he will be making. Basically, starting from a position with value $\omega^2$, black plays to a position with value at least $\omega\cdot n$ for some $n$. A position with value at least $\omega^2\cdot m$ allows black to make $m$ many such large announcements, where each large announcement is the number of announcements to be made, with each announcement being a lower bound on the number of moves until the next announcement, and the next large announcement occurring when all such smaller announcements are completed. In a position with value at least $\omega^3$, black starts with a huge announcement, which is the number $m$ of large announcements to be made, by playing to a position with value at least $\omega^2\cdot m$. A position with value $\omega^3\cdot r$ allows black to make $r$ many such huge announcements, and a position with value $\omega^4$ allows black to announce at the start an arbitrarily large finite number $r$, his meta-announcement, that will be the number of huge announcements.

Perhaps it will be helpful to illustrate how this announcement manner of understanding the game values works with specific game values. In effect, black is counting down in the ordinals, since white always plays to reduce the game value by exactly one, and black must choose a smaller ordinal at every limit stage. Starting from $\omega^4$, black might play to a position with value, say, $\omega^3\cdot 5+\omega^2\cdot 17+\omega\cdot34+1234$. This amounts to making a meta-announcement of $6=5+1$, with the first huge announcement being $18=17+1$, and the first large announcement being $35=34+1$, and the first announcement being $1234$. The game values proceed to count down in $1234$ moves to $\omega^3\cdot 5+\omega^2\cdot 17+\omega\cdot 34$, at which time black makes his second announcement by playing to $\omega^3\cdot 5+\omega^2\cdot 17+\omega\cdot 33+2^{10}$, say, which announces another $2^{10}$ moves before the next announcement. Proceeding in this way, black will make $33$ such additional announcements, with corresponding play after each of them, before arriving at the ordinal $\omega^3\cdot 5+\omega^2\cdot 17$ itself, when it is time to make another large announcement. At this point, let us suppose black plays to a position with value $\omega^3\cdot 5+\omega^2\cdot 16+\omega\cdot 10^{100}+2^{2^{10}}$, which amounts to a large announcement of $10^{100}+1$ and a first announcement in that cycle of $2^{2^{10}}$. And so on. Eventually, black will use up the $10^{100}$ announcements in this cycle, and then ultimately all $16$ additional large announcements in this huge cycle, arriving at $\omega^3\cdot 5$, at which time he uses another of his precious five remaining huge announcements by playing to $\omega^3\cdot 4+\omega^2\cdot 10^{100!}+\omega\cdot 10^{10^{100}}+2^{2^{2^{100}}}$, say, which makes a huge announcement of $10^{100!}+1$, with the first large announcement of $10^{10^{100}}+1$ in that cycle and the corresponding first announcement of $2^{2^{2^{100}}}$ in {\it that} cycle. Eventually, all $5$ huge announcements will be used up and all corresponding large announcements, and all announcements, and in the very last leg of play, black will be counting down from some very large natural number one by one until $0$ is reached, at which point white wins.

\section{A position with value $\omega^4$}\label{Section.APositionWithValueOmega^4}

Let us now describe our new position, shown in figure~\ref{Figure:MainPositionValueOmega^4}, which has
value $\omega^4$, with black to play. The position has three principal components: the rook tower block
at the upper right, which is adapted from the Evans-Hamkins position of
\cite{EvansHamkins2014:TransfiniteGameValuesInInfiniteChess}; the throne room at the upper left, where
the two kings sit facing each other in uneasy if ultimately temporary d\'etente; and the bishop cannon
battery arrangement on the bottom half of the board, consisting of the bishop cannons at the lower right
and the gateway terminals at the lower left, with a large open arena between them. Although the figure
shows only finitely much of the position, we intend for the position we are describing that the indicated patterns should continue indefinitely. In particular, the position has infinitely many pieces; the pattern of
rook towers continues infinitely toward the upper right; the bishop cannons appear one after another,
with one cannon of each finite size, continuing toward the lower right; and the the gateway terminals
continue with one wing of each finite size toward the lower left.

In large scale, during optimal play, what will happen is that white makes slow but steady progress in the rook tower arrangement at the upper right, just as in the original Evans-Hamkins position, with the focus of play moving from an initial active tower, determined by the very first move of the key bishop up on the diagonal, and then continuing successively to the adjacent towers to the left, until finally a bishop is released to the left of the left-most rook tower, which will deliver checkmate to the black king in the throne room; at nearly every step of this progression, however, black uses the bishop cannon battery to mount increasingly desperate distractions, prolonged bouts of doomed yet forced play, in which successive series of black threats must be answered by white. Thus, the bishop cannon battery replaces the simple black rook white king harassment arrangement in the original Evans-Hamkins position~\cite[fig~10]{EvansHamkins2014:TransfiniteGameValuesInInfiniteChess}, while giving black a higher order of delay.

Specifically, in the original Evans-Hamkins position~\cite[fig~10]{EvansHamkins2014:TransfiniteGameValuesInInfiniteChess}, there were two components, the main block of rook-towers, which by itself had a game value of $\omega^2$, and the rook-king harassment arena, which enabled black to delay progress in the main line with a multiplier of $\omega$, leading to an overall game value of $\omega^3=\omega\cdot\omega^2$ for that position. In the new position here, we essentially replace the rook-king harassment arena with the bishop cannon battery, which enables black repeatedly to delay white's progress in the main line with a multiplier of $\omega^2$, rather than merely $\omega$ as before, and this leads to an overall game value of $\omega^4=\omega^2\cdot\omega^2$ for the new position.

Let us now analyze the position in closer detail, considering each of the main components in turn. After this, we shall describe the main line of optimal play, along with the ordinal analysis, and explain why deviations from the main line of play are disadvantageous.

\begin{figure}[h]
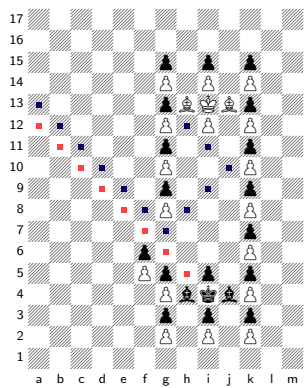

\catcode`\.=\active\def.{{1}}
\chessboard[maxfield=m17,
            boardfontsize=8pt,
            label=false,
            labelbottom=true,
            labelleft=true,
            labelfontsize=5pt,
            labelleftwidth=2ex,
            hlabelformat={\hfill\arabic{ranklabel}\,\,},
            showmover=false,
            border=false,
            margin=false,
            pgfstyle=color,
            color=red!75!white,
            padding=-3pt,
            backfields={a12,b11,c10,d9,e8,f7,g6,h5},
            pgfstyle=color,
            color=blue!50!black,
            padding=-3pt,
            backfields={a13,b12,c11,d10,e9,f8,g7,h8,i9,j10,i11,h12},
            setfen=%
.........../%
.........../%
......p.p.p.../%
......P.P.P.../%
......pBKBp.../%
......P.P.P.../%
......p...p.../%
......P...P/%
......p...p./%
......P...P/%
..........p./%
.....p....P/%
.....Pp.p.p./%
......PbkbP/%
......p.p.p/%
......P.P.P/%
........../%
]%
\caption{The throne room}\label{Figure:ThroneRoom}
\end{figure}

Consider first the throne room, pictured in the detail of figure~\ref{Figure:ThroneRoom}. A close
inspection of this arrangement will reveal that it is completely locked up; none of the pieces can
legally move. Meanwhile, both of the kings are exposed to checkmate, if a bishop of the opposite color
should position itself correctly. In particular, a white bishop can deliver checkmate to the black king
from any square on the diagonal marked in red, such as {\tt 1.Bc10\#}.  In contrast, a black bishop can
deliver checkmate to the white king by following the the indicated blue zig-zags via {\tt 1\ldots Bc11\
2\ldots Bg7\ 3\ldots Bj10\ 4\ldots Bh12\#}. Any white-square white bishop or black-square black bishop
that is free in the open zone around the throne room can deliver checkmate in just a few moves, and this
checkmate can be delivered by means of the indicated diagonals leading to the door of the throne room.
Notice the key feature of the position, namely, that while the white bishop checkmates directly from that
diagonal, the black checkmate takes exactly three additional steps after arriving on that diagonal. This
precise difference will be important in our analysis of the threats made by black in the gateway terminal
configuration.

\begin{figure}[h]
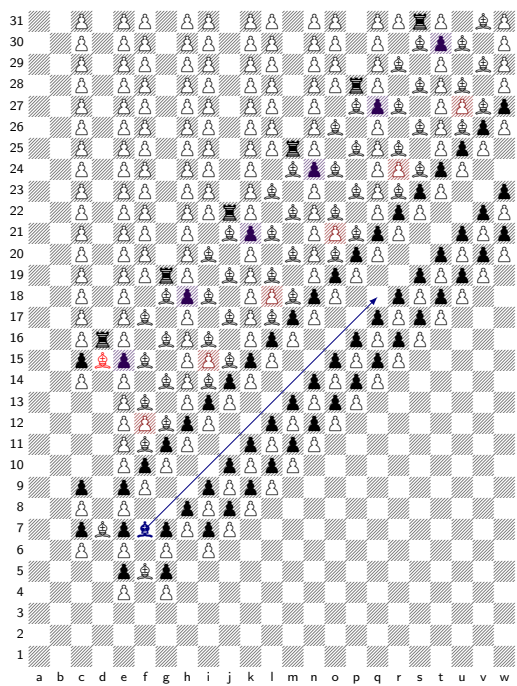

\catcode`\.=\active\def.{{1}}
\chessboard[maxfield=w31,
            boardfontsize=8pt,
            label=false,
            labelbottom=true,
            labelleft=true,
            labelfontsize=5pt,
            labelleftwidth=2ex,
            hlabelformat={\hfill\arabic{ranklabel}\,\,},
            showmover=false,
            border=false,
            margin=false,
            coloremphstyle=\color{blue!65!red!50!black},
            emphfields={e15,h18,k21,n24,q27,t30},
            coloremphstyle=\color{red!50!black},
            emphfields={f12,i15,l18,o21,r24,u27},
            coloremphstyle=\color{blue!50!black},
            emphfield=f7,
            pgfstyle=straightmove,
            linewidth=10000sp,
            shortenstart=-2pt,
            color=blue!50!black,
            backmove={f7-q18},
            coloremphstyle=\color{red},
            emphfield=d15,
            setfen=%
..P.PP.PP.PP.PP.PPrP.BP/%
..P.PP.PP.PP.PP.P.BpB.P/%
..P.PP.PP.PP.PP.PB.P.BP/%
..P.PP.PP.PP.PPrP.BPB.P/%
..P.PP.PP.PP.P.BpB.PPBp/%
..P.PP.PP.PP.PB.P.BPBpP/%
..P.PP.PP.PPrP.BPB.PpP./%
..P.PP.PP.P.BpB.PPBpP../%
..P.PP.PP.PB.P.BPBpP..p/%
..P.PP.PPrP.BPB.PpP..pP/%
..P.PP.P.BpB.PPBpP..pPp/%
..P.PP.PB.P.BPBpP..pPpP/%
..P.PPrP.BPB.PpP..pPpP/%
..P.P.BpB.PPBpP..pPpP./%
..P.PB.P.BPBpP..pPpP../%
..PrP.BPB.PpP..pPpP.../%
..pBpB.PPBpP..pPpP/%
..P.P.BPBpP..pPpP./%
....PB.PpP..pPpP.../%
....PPBpP..pPpP..../%
....PBpP..pPpP..../%
....PpP..pPpP............/%
..p.pP..pPpP............../%
..P.P..pPpP............/%
..pBpbpPpP.........../%
..P.P.P.P............/%
....pBp................/%
....P.P.........../%
.........../%
.........../%
]%
\caption{The rook towers}\label{Figure:RookTowers}
\end{figure}

Next, consider the rook tower configuration, shown in the detail of figure~\ref{Figure:RookTowers}. This
part of the position is drawn directly from the
corresponding part of the Evans-Hamkins position of
\cite{EvansHamkins2014:TransfiniteGameValuesInInfiniteChess}, with some changes along the left and bottom edges, and we refer the reader to that article for additional thorough discussion. Notice that the black bishop at {\tt f7} (highlighted in blue) is attacked by the white pawns at {\tt e6} and {\tt g6}, and if it should be captured by one of those pawns,
then one of the white bishops at {\tt d7} or {\tt f5} will be released into the free zone, where it can
deliver checkmate very quickly in the throne room. The reader may similarly verify that if black were to
capture any of the pawns at {\tt e6, e8} or {\tt g6} with the {\tt f7} bishop, then white will recapture
with a bishop and again be free to checkmate quickly. Since black has no other such immediate threats
anywhere on the board, the main line of play will be that on the very first move, black will move his
bishop at {\tt f7} up the diagonal as far as he wishes. In our ordinal analysis, we will argue that for
black to move the bishop out past $n$ rook towers will result in a position with value between
$\omega^3\cdot n$ and $\omega^3\cdot(n+1)$, and it is to black's advantage (in terms of maximizing the
game value) to move the {\tt f7} bishop out as far as he can (but alas, he must move only some specific
amount); and this is why the overall position, before the {\tt f7} bishop move, has value $\omega^4$.

After the black bishop moves up the diagonal, the bishop will immediately be captured by a white pawn, such as
{\tt 1\ldots Bq18 2.rxq18}. This white pawn will advance to pierce the upper pawn diagonal protecting the
rook towers, with {\tt 3.q19, 4.qxp20, 5.pxq21, 6.p20}, aiming to release the bishop chain
starting with {\tt 7.Bq18}. Note that black gains no advantage by capturing the white pawn with {\tt 3\ldots pxq19}, since this will only release the bishops more quickly. (Instead, black will use all the free black moves in this sequence to mount long bouts of forced-reply moves using the bishop cannons, as we shall explain.) After a few more
bishop moves, white will be able to advance the key pawn at {\tt o21} (in red), which will advance and
attack the black guard pawn at {\tt n24} (in purple). Black will not take the white attacking pawn, since
this will simply release the bishops in the next tower (to the left) sooner than necessary. Instead, the
full line of play is {\tt 1\ldots Bq18 2.rxq18, 3.q19, 4.qxp20, 5.pxq21, 6.p20, 7.Bq18, 8.Bpo20, 9.B22p21, 10.o22, 11.o23, 12.oxn24 Rm{\it j} 13.lxm{\it j}}, so that the black rook at {\tt m25}
is attacked and then moves up the rook tower to an arbitrary height $j$, the higher the better, and then
is taken from the left by a white pawn. This opens up a hole in the pawn column in file {\tt l}. White will
proceed to move those pawns up one by one until after a very long time, the bishop at {\tt m24} is able
to move into the hole created at {\tt l25}. After this, the bishops in the {\tt l} and {\tt m} files can
advance upward in such a way that the key pawn at {\tt l18} is able to advance and attack the next rook
tower guard pawn at {\tt k21}, thereby activating the next rook tower. This whole process repeats for
each rook tower down the line to the left: bishops are enabled to move out of the way of the key white
pawns (red), which advance to attack and take the black guard pawns (purple), thus attacking the black
rook of the next tower; the black rook advances very high, is captured from the left, opening a hole in
the pawn column; white advances pawns until the bishops of the next tower are enabled to release the next
key pawn. Finally, the left-most rook tower is activated when the white key pawn at {\tt f12} advances
(since the bishops will have moved out of the way) to attack and take the black guard pawn at {\tt e15},
causing the last rook at {\tt d16} to advance to arbitrary height. White takes from the left, creating a
hole in the pawn column on the {\tt c} file. White advances these pawns one by one until finally the pawn
at {\tt c16} advances, which releases the mating bishop at {\tt d15} (blood-red) into the free zone; this
bishop
will quickly deliver checkmate to the black king in the throne room.

Thus, white is able to make systematic progress in the rook tower, which will lead inevitably to a white bishop being released to deliver checkmate. A key feature of the main line of play in the rook tower is that black has numerous free moves. For example, whenever one of the black rooks advances and is captured in the main line of play, white subsequently advances an enormous number of pawns in that column in order to move the hole in the pawn column down, so as to release the bishops in the adjacent rook tower to the left, thereby activating the next tower. All the while that white is advancing these pawns, black is free to play elsewhere. What black will do at any such opportunity is initiate firing with another (arbitrarily large) bishop cannon, which will lead to an extremely long bout of forced-play in the bishop cannon battery arena, a doomed but time-consuming distraction for white. Indeed, in optimal play, every single advance of white in the rook tower comes at the price of yet another dramatic spectacle mounted by black in the bishop cannon battery.

\begin{figure}[h]
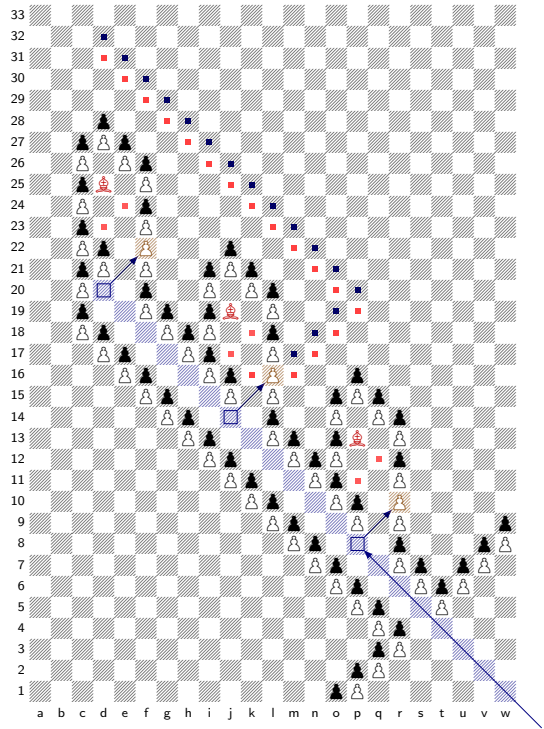

\catcode`\.=\active\def.{{1}}
\chessboard[maxfield=w33,
            boardfontsize=8pt,
            label=false,
            labelbottom=true,
            labelleft=true,
            labelfontsize=5pt,
            labelleftwidth=2ex,
            hlabelformat={\hfill\arabic{ranklabel}\,\,},
            showmover=false,
            border=false,
            margin=false,
            coloremphstyle=\color{blue!50!black},
            emphfields={d20,e19,f18,g17,h16,i15,j14,k13,l12,m11,n10,o9,p8,q7,r6,s5,t4,u3,v2,w1},
            coloremphstyle=\color{red!75!black},
            emphfields={d25,j19,p13},
            coloremphstyle=\color{orange!50!black},
            emphfields={f22,l16,r10},
            pgfstyle=color,
            color=red!65!white,
            padding=-200000sp,
            backfields={e24,d23,k18,j17,q12,p11},
            pgfstyle=border,
            padding=-100000sp,
            linewidth=20000sp,
            color=blue!50!black,
            backfields={d20,j14,p8},
            pgfstyle=straightmove,
            linewidth=20000sp,
            shortenstart=-20pt,
            shortenend=3pt,
            color=blue!50!black,
            backmove={w1-p8},
            pgfstyle=straightmove,
            linewidth=20000sp,
            shortenstart=3pt,
            shortenend=4pt,
            color=blue!50!black,
            backmove={d20-f22,j14-l16,p8-r10},
            pgfstyle=color,
            color=blue!50!black,
            padding=-3pt,
            backfields={m17,n18,o19,p20,o21,n22,m23,l24,k25,j26,i27,h28,g29,f30,e31,d32},
            pgfstyle=color,
            color=red!75!white,
            padding=-3pt,
            backfields={k16,m16,n17,o18,p19,o20,n21,m22,l23,k24,j25,i26,h27,g28,f29,e30,d31},
            setfen=%
............./%
............./%
............./%
............./%
............./%
...p........../%
..pPp.........../%
..P.Pp........../%
..pB.P........../%
..P..p........../%
..p..P......./%
..Pp.P...p.../%
..pP.P..pPp..../%
..P..p..P.Pp.../%
..p..Pp.pB.P..../%
..Pp..PpP..p.../%
...Pp..Pp..P..../%
....Pp..Pp.P...p...../%
.....Pp..P.P..pPp../%
......Pp...p..P.Pp./%
.......Pp..Pp.pB.P...../%
........Pp..PpP..p..../%
.........Pp..Pp..P......p....../%
..........Pp..Pp.P.....pP...../%
...........Pp..P.P....pP..../%
............Pp...p...pP...../%
.............Pp..Pp.pP...../%
..............Pp..PpP..../%
...............Pp..P...../%
................Pp....../%
................pP...../%
...............pP...../%
..............pP...../%
.............pP...../%
............pP...../%
]%
\caption{The bishop gateway terminal wing of length three}\label{Figure:BishopGateways}
\end{figure}

\begin{figure}[h]
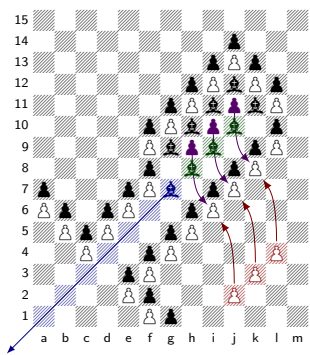

\catcode`\.=\active\def.{{1}}
\chessboard[maxfield=m15,
            boardfontsize=8pt,
            label=false,
            labelbottom=true,
            labelleft=true,
            labelfontsize=5pt,
            labelleftwidth=2ex,
            hlabelformat={\hfill\arabic{ranklabel}\,\,},
            showmover=false,
            border=false,
            margin=false,
            coloremphstyle=\color{blue!50!black},
            emphfields={a1,b2,c3,d4,e5,f6},
            pgfstyle=straightmove,
            linewidth=10000sp,
            shortenstart=-2pt,
            shortenend=-20pt,
            color=blue!50!black,
            backmove={g7-a1},
            coloremphstyle=\color{green!25!black},
            emphfields={h8,i9,j10},
            coloremphstyle=\color{blue!50!black},
            emphfield=g7,
            coloremphstyle=\color{red!50!black},
            emphfields={j2,k3,l4},
            coloremphstyle=\color{blue!55!red!75!black},
            emphfields={h9,i10,j11},
            pgfstyle={[clockwise=false,style=knight]curvemove},
            linewidth=20000sp,
            shortenstart=0pt,
            shortenend=3pt,
            color=blue!55!red!75!black,
            backmove={h9-i6,i10-j7,j11-k8},
            pgfstyle={[clockwise=false,style=knight]curvemove},
            linewidth=20000sp,
            shortenstart=2pt,
            shortenend=5pt,
            color=red!50!black,
            backmove={j2-i6,k3-j7,l4-k8},
            setfen=%
.............../%
.........p..../%
........pPp..../%
.......pPbPp/%
......pPbpbP/%
.....pPbpb.p../%
.....Pbpb.pP../%
.....p.b.pP.../%
p...pPb.pP../%
Pp.pP..pP....../%
.PpP..pP...../%
..P..pP....P/%
....pP....P/%
....Pp...P/%
.....Pp../%
......Pp..../%
]%
\caption{A three-shooter bishop cannon}\label{Figure:BishopCannon}
\end{figure}

So let us now describe how play proceeds in the bishop cannon battery arrangement, which consists of the
bishop cannons at the lower right and the gateway terminals at the lower left, separated by the large
enclosed arena between them. A wing of bishop gateway terminals is detailed in figure
\ref{Figure:BishopGateways} and a bishop cannon in figure~\ref{Figure:BishopCannon}. The entire battery
configuration is sealed-in by the pawn chain surrounding the lower part of the board, and the pattern of
play here is that black will make various threats for a black bishop to escape into the free area
containing the throne room, where it would be able to make checkmate. But in fact, white will be able to
answer all such threats, and so in fact no bishop will actually escape. More specifically, at any stage
of play, black will be able to initiate firing from one of the bishop cannons, of any desired size. The
initial firing of any cannon is not itself an immediate threat, and so immediately after firing a cannon
in such a way for the first time, white will be able to make progress in the main line, such as by
advancing a pawn in the currently active rook tower; and indeed, such kind of isolated moves will be
nearly the only opportunity that white will have for making progress on the main line (which itself
already takes some time for white, with game value $\omega^2$ even apart from the bishop cannons). Once a
cannon is fired for the first time, then black will subsequently be able to fire it repeatedly, each time with force, so that white must respond, until its ammunition is depleted; each time the cannon is fired, a black bishop moves out the exit diagonal of the cannon, stopping at any desired arbitrarily large entrance diagonal to one of the
bishop gateway terminals, and white responds on the cannon side. By next entering the gateway terminal,
black is able to make a series of successive threatening moves inside that particular gateway terminal
wing. When all such threats are used up, the next bishop fires out of the cannon, which white must answer
on the cannon side, after which the bishop makes a series of threats inside whichever gateway terminal
that bishop had chosen. And so on, until all the bishop ammunition from that cannon is depleted. At this
stage, black may chose to initiate firing with another (much larger) cannon, and only then does white get the one-move
reprieve for progress in the main line.

Let us look closely at the line of play inside the gateway terminals. Consider the bishop gateway terminal wing of length three pictured in figure~\ref{Figure:BishopGateways}. The position looks something like an airport terminal wing, with a long diagonal hallway (blue), having several gates (in this case three gates) leading off from it. The hallway has several attacking squares for the black bishop (marked with blue squares), and each attacking square attacks a gateway door (brown), and each gateway has a guard tower with a bishop guard (red). We are to imagine that as a result of one of the bishop cannons firing, a black bishop has placed itself out on the blue entrance diagonal to this particular gateway wing. We claim that for black to move his bishop to any of the three attacking squares {\tt p8, j14}, or {\tt d20} is an immediate threat that white must answer. To see this, let us suppose that black has placed a bishop on {\tt j14} with {\tt 1\ldots Bj14}. The correct and forced reply, we claim, is for white to summon the bishop guard with {\tt 2.Bk18\forcedmove}, where we use the symbol \forcedmove\ to indicate that this move is forced. If black now takes the gateway pawn with {\tt 2\ldots Bxl16}, then white follows with {\tt 3.Bj17\forcedmove}. If black should actually exit the gateway, following the pathway marked in blue, then the white guard bishop will follow closely, two steps behind, on the red pathway with {\tt 2\ldots Bp20 3.pl16 Bd32 4.Bl15, 5.Bp19}. The point is that as the black bishop races to the throne room, the white guard bishop can always land on a corresponding square directly underneath, but two moves later. In particular, when black arrives on the diagonal entryway to the throne room (discussed in connection with the throne room above), white will arrive on that diagonal two moves later. But since, as we discussed earlier, black requires three moves after landing on that diagonal in order to deliver checkmate, white's mere two-move delay is sufficient, and white will deliver checkmate first. Thus, black should not actually exit the gateway, if the white guard had moved when the gate was first attacked. But if white had ignored the original threat posed by black's move {\tt 1\ldots Bj14}, then white would be three moves behind, rather than only two, and in this case black would be able to checkmate first. Because of this, for black to place a bishop on one of the attacking squares in the gateway terminal is indeed a threat that must be answered immediately by white summoning the bishop guard, which moves one step closer to the gateway.

The main point of the arrangement is that black is therefore able to make a series of such threats
successively, simply by moving from one attack square to the next along the terminal diagonal. To
illustrate, starting from the original position of figure~\ref{Figure:BishopGateways}, with a black
bishop out on the blue diagonal, then optimal play would be {\tt 1\ldots Bp8 2.Bq12\forcedmove\ Bj14 3.Bk18\forcedmove\ Bd20 4.Be24\forcedmove\ Bxf22 5.Bd23\forcedmove}, and now black moves on to other projects, making sure never to move his bishop on {\tt f22}. Note that black could have used the attack squares {\tt p8}, {\tt j14} or {\tt d20} in a different order, although each such attack square can be used at most once. Indeed, there are precisely six identical optimal lines for black, each following this same form and each taking precisely five moves. On the terminal wing having one hundred gates, black will be able to make at least one hundred forced-reply moves in succession, plus one more by occupying the doorway, as we have indicated; with $n$ gates, black can make $n+1$ such forced-reply moves.
To summarize how play proceeds in a gateway, the main point is that for black to have a bishop on the
entrance diagonal of a gateway terminal wing of length $n$ will enable black to make $n+1$ successive
immediate threats, each of which must be answered by white in a way that does not affect the relevant
operation of the other parts of the position.

Consider now the firing of the cannons. Figure~\ref{Figure:BishopCannon} details a three-shooter cannon,
which will mean that after the initial firing of the cannon, which is not a forced-reply move, black will be
able to fire three additional bishops out of the cannon with force, so that white must reply immediately
on the cannon side (so four shots in all, but only three with forced-reply). To initiate firing with this
particular cannon, black will move the front bishop at {\tt g7} out on the blue diagonal, placing it on a
square that faces an as-yet unused gateway terminal entrance diagonal, so that it will subsequently be
able to enter that gateway terminal. For black to fire the cannon by moving the front bishop out like
this is not an immediate threat to white, who may at this moment play elsewhere in order to make progress
on the main line, such as by advancing a pawn in the currently active rook tower. Indeed, in the full
line of optimal play, the only chance that white gets to make progress on the main line in the rook tower
configuration is in such a single move after black first fires a new cannon (except in the last four or
five moves of the full line, when black's threats are not fast enough). After the cannon is thus
initiated and white has replied in the rook tower, then black may use his front bishop (now out in the
arena facing a gateway terminal) to make a series of forced-reply threats inside that terminal. After
this, he may fire the cannon again by moving the next bishop at {\tt h8} out of the cannon onto the blue
diagonal, making sure to position it on the entrance diagonal of an as-yet-unused gateway terminal, the
larger the better. The key feature of the cannon is that all of the bishops fired after the front bishop
move out of the cannon with force, making an immediate threat to which white must reply. What we claim
for the current cannon is that immediately after the bishop at {\tt h8} is fired, white must advance the
guard pawn on {\tt j2}, aiming to protect the white pawn on {\tt i6}. If white does not advance that
guard pawn, then
black will advance his pawn at {\tt h9}, and in one additional move will be attacking the (insufficiently defended)
pawn on {\tt i6}. If black captures the pawn on {\tt i6} and is not immediately retaken, then black can
advance the capturing pawn and then advance the pawn on {\tt i7}, which will open a hole in the pawn
chain, allowing black bishops to escape into the free area, where they will quickly deliver checkmate in
the throne room. Note that it does no good for white to capture the attacking pawn on {\tt h7}, since in
this case black can simply advance the {\tt i7} pawn, again opening a hole. Let us consider the critical line: {\tt 1\ldots B$xy$ {\rm (from {\tt h8})} 2.j3\forcedmove\ h8 3.j4\forcedmove\ h7 4.j5\forcedmove\ hxi6 5.jxi6\forcedmove}, which guards against the threat {\tt 5\ldots i5, 6\ldots i6}, winning by opening lines for the bishop on {\tt g9}. To summarize, after moving the bishop at {\tt h8} out onto the diagonal, the forced play can proceed with {\tt 1.Pj3 Ph8 2.Pj4 Ph7 3.Pj5 Pxi6 4.Pxi6}, which shores up the threat that had been posed by the discovered threat of {\tt h9}
after the bishop on {\tt h8} had fired out of the cannon. After this sequence, black is now free to use
his bishop (which is now facing a fresh gateway terminal) to make an additional series of threats inside
that gateway terminal. After that, he may fire the next bishop out of the cannon, moving from {\tt i9}
out on the diagonal to face a new (and much larger) gateway terminal, the larger the better. This again
is a threat, on account of the pawn at {\tt i10}, and so black will advance the pawn at {\tt k3} to
defend. Thus, each of the bishops marked in green in the cannon fires with force, in the sense that white
must reply by advancing the corresponding protection pawn (red) in order to shore up the boundary wall.
Black will have no incentive to fire the bishop at {\tt k11} out of the cannon, since this is a move to
which white need not reply immediately; and so black is better off opening up a new cannon than firing
that bishop. In summary, when black initiates a new cannon of size $b$ by firing it for the first time,
moving its front bishop out, then this initial firing is not itself an immediate threat, but after white
responds elsewhere, black will be in a position of being able to fire $b$ bishops out of the cannon, each
with forced-reply, and each of these bishops will be able to attack any desired unused gateway terminal
of arbitrarily large size $g$, after which it will be able to make $g$ additional forced-reply threats
inside the terminal. Thus, after the initial firing of a cannon of size $b$, black is enabled to
undertake forced-reply play of size $\omega\cdot b$. Furthermore, a key feature of the entire bishop
cannon battery is that it is {\it reusable}, in the sense that once black has initiated fire with a
bishop cannon, eventually depleting it, then at any opportunity black may choose to initiate firing with
another cannon, of arbitrarily large size, the larger the better. The fact that some cannons and gateway
terminals (only finitely many) have already been used does not compromise in any way the future use of
the other much-larger cannons and gateway terminals.

Let us now summarize what we claim is the main line of best play for the overall position. The position calls for black to play first. Black will move the key bishop in the rook towers up on the diagonal as far as he cares to move, the farther the better. White captures that bishop with a pawn, and aims to advance that pawn in order to disrupt the inner diagonal wall, releasing the bishops of the rook tower in that vicinity. White then aims to advance the key pawn that had been sealed-in by those bishops, in order to attack and capture the black guard pawn and ultimately attack the black rook of that tower, which black will move up as high as he cares to move, where it will be captured by a white pawn, creating a hole in the pawn column, which white aims to steadily move downward, in order to release the bishops of the adjacent rook tower to the left, where the whole process starts again with that tower. At nearly every step of play, however, whenever given the opportunity, black initiates firing with a fresh cannon in the bishop cannon battery, as large as he cares, the larger the better, and proceeds to make a long series of forced-reply moves using the cannon and gateway terminals, in the manner described in the previous paragraph. The only chance white has to make progress in the rook tower configuration, therefore, is immediately after black moves so as to initiate firing with a fresh cannon (also, there are a few white moves in the rook tower where black must reply there, such as on the very first move, when he must move his bishop, and also in each rook tower, when the rook is finally attacked, black must move it up). Eventually, when the currently active rook tower has moved all the way to the left-most rook tower and when that rook tower is completed, then the key bishop at the bottom left of that tower is released, which will proceed to deliver checkmate to the black king in the throne room. At the very end of play in the final tower, once the bishop is nearly released, then the bishop cannon threats are no longer operative, since they are not quick enough, and so in the very final stages of play, black's threats will go unanswered, while white makes the final four or five moves to checkmate.

Let us argue next that this pattern of play means that the ordinal game value of the position is precisely $\omega^4$. We will subsequently use this ordinal analysis to explain why deviations from the main line of play are suboptimal. Suppose that we are in the midst of carrying out the main line of play as we have just described it (but let us suppose that we have not yet reached the final four-or-five moves where the bishop cannon threats no longer carry force). Let us associate to the current position the ordinal $$\omega^3\cdot r+\omega^2\cdot h+\omega\cdot b+ g,$$ where $r$ is the number of remaining rook towers not yet activated, $h$ is the number of moves left for white in the current rook tower (counting all the moves up until the rook of the next tower is threatened), $b$ is the number of bishops that we have left in the currently active cannon (counting only those bishops that are capped by a pawn, as in green in figure~\ref{Figure:BishopCannon}, so that they fire with forced-reply), and $g$ is the number of forced-reply threats available from the current bishop in the active gateway terminal (this is the number of gateways left in the currently active gateway terminal, plus a few extra moves because of the threats we discussed that black can make by actually
entering the last gateway and by advancing the discovered pawn inside the cannon). Notice that as black
uses up the threats available from the current bishop in the active gateway terminal, the number $g$
steadily descends by one each move. And when $g=0$, then black fires the next available bishop out of the
cannon, which reduces $b$ by one, but now allows $g$ to become as large as desired, since black can move
that bishop out to attack any arbitrarily large unused gateway terminal. When the current cannon is out
of ammunition (counting only those bishops that fire with forced-reply), then black can initiate fire
with a new cannon, which allows white to decrease $h$ by one, since white makes progress in the rook
tower when black engages the cannon, but then $b$ can become arbitrarily large, since black can engage a
cannon of arbitrarily large size. When $h$ becomes zero, then the next rook tower will become activated,
which decreases $r$ by one, but allows the new $h$ value to be arbitrarily large, since black can move
his rook up as high as desired. In summary, the main line of play allows black to count down in the
ordinals as we have described, and so provided that neither player has incentive to deviate from the main line, it follows that the ordinals that we have associated with the positions in the main line are precisely the ordinal values of those positions, leading to value $\omega^4$ overall for the full position.

Naturally, we must argue that neither player gains advantage by deviating from this main line of play, and that in particular, white cannot force a quicker checkmate---quicker in the sense of having a smaller ordinal game value---and black has no winning strategy and furthermore cannot impose a larger ordinal game value for white than we have described. To show this, we shall consider various kinds of unexpected play and argue that neither player has
an interest in undertaking it.

Let us first consider the various possibilities for black to pursue alternative lines of play and argue that it is not in black's interest to do so. We claim in particular that black cannot force checkmate and furthermore has no line of play realizing a higher ordinal game value than $\omega^4$.

To begin, let us consider that black may deviate from the main line by amassing an army of free black-square bishops in the bishop cannon arena, simply by moving the initial bishops out of more and more cannons, without using them directly as threats in the bishop gateways. Black may amass as many bonus bishops in this way as he wishes, simply by using up one rook tower in order to free up enough time to move the bishops out. We could even imagine allowing black freely to place extra black-square black bishops in the arena area, as long as they are not placed directly onto the attacking squares in the bishop gateways. We claim, nevertheless, that such a bishop army is powerless to exert any extra threat on white, and black is better off following the main line. The reason is that the black-square black bishops are penned in by the boundary wall of black-square black pawns that enclose the bishop-cannon arena. The only squares that black-square black bishops can attack in the bishop cannon arena are the gateway doorway pawns discussed in figure~\ref{Figure:BishopGateways}, and white need reply to black's moves in the cannon arena only when a bishop actually attacks such a gateway doorway pawn, or when one of the bishop cannons fires with force (by moving a pawn-capped bishop out of a cannon) as described in figure~\ref{Figure:BishopCannon}. White has nothing to fear from extra free black-square black bishops in the arena, until they attack a doorway pawn, in which case the threat can be immediately answered by white by moving the corresponding white bishop guard into place, as discussed in connection with figure~\ref{Figure:BishopGateways}.

Next, let us consider what might be a worrisome feature for white in the main line of play, namely, the feature that black bishops will eventually occupy the doorways of numerous bishop gateways. Our analysis in connection with figure~\ref{Figure:BishopGateways} ends with each black bishop actually taking a gateway doorway pawn and occupying the doorway, but not exiting the gateway, with the white guard bishop in place, ready to pursue (two steps behind), should the bishop actually exit the doorway. Is there any way that black can make use of having an enormous number of such almost-escaped bishops in waiting? No. It doesn't matter how many black bishops are ready to exit, because once the first one leaves its gateway doorway, then the main line calls for white to advance the white pawn below the doorway and thereby allow the white guard bishop to escape, with fewer moves to checkmate than black. No move by black to use another doorway-occupying bishop would prevent white from delivering checkmate. Since the white bishops move exclusively on white squares and black on black squares, there is no means for black to place obstacles in white's checkmate path using extra black-square black bishops. So there is no danger to white in black having a large number of almost-escaped black bishops occupying the doorways, provided that the corresponding white guard bishops are in place for each one of them.

Next, inside one of the bishop cannons itself, let us consider that black attempts an alternative line of play with his black pawns, after the bishops are fired out of the cannon. Black would be free to fire the bishops out of the cannon, without actually moving the pawns above in the way that was described in connection with figure~\ref{Figure:BishopCannon}. In this way, he would have all the bishops cleared out of the way, before moving any of the pawns. Is there any advantage in this? No. If, as called for in the main line, white has moved his guard pawns to protect the wall, then it doesn't matter if black tries to coordinate his cannon pawns. The only way for black to make use of the pawns is to perform a capture of the corresponding wall pawn, which will then be immediately taken back by the white guard pawn, reinforcing the wall, and there is no way for the black pawns to reinforce each other in this process or to gain the support of a black-square black bishop. Thus, black is better off using each pawn before firing the next bishop, since this preserves more ordinal value by reducing only $g$ instead of $b$.

Let us now consider the possibility that black tries to use the bishop cannons in a different order than we have described in the main line. Can he somehow delay his fate for longer? To analyze this situation, we refine the definitions of the parameters $b$ and $g$ in a way that is consistent with their definitions in the main line but accommodates certain deviations by black. Specifically, let $b$ be the total number of pawn-capped bishops left in \emph{all} currently active cannons (in case black has more than one cannon active at once), and let $g$ be the maximum number of consecutive moves that black can make (possibly using many different pieces), each of which either threatens to
capture a gateway pawn or actually captures a gateway pawn or moves an unblocked pawn so that white must
defend with a guard pawn. Once black has initiated fire with one of the bishop cannons, it is in black's best interest to continue making as many forced-reply moves as possible, reducing $g$ by one each time, or, if $g=0$, reducing $b$ by one and increasing $g$ by a large number, since otherwise white will take the opportunity to reduce $h$ by one, by making progress in the rook tower, and this is a larger reduction in the ordinal than black could have ensured by making a forced-reply move to reduce only $g$ or $b$. Similarly, whenever a bishop cannon is depleted, then it is in black's interest to initiate fire with another cannon, the larger the better, since this allows him to increase $b$ by a large amount, whereas any other bishop move (such as bringing a non-pawn-capped bishop from a depleted gun into firing position or retreating a bishop out of a terminal for redeployment) would at best increase $b$ by one (or threaten to do so), while giving white an opportunity to decrease $h$ by one, which is a far greater reduction in the ordinal and therefore worse for black.

Let us next consider the possibility that black attempts to amass an army of rooks inside the rook tower. For example, black can simply move one of the black rooks in a rook tower to capture an adjacent white pawn or the white bishop below it. Although white can simply capture the rook immediately, this situation is actually more dangerous for white than it may initially appear. The reason is that even though white can immediately answer all such crazy rook moves, nevertheless white cannot afford to adopt the strategy of {\it always} answering such moves, since then black could force a draw via infinite play, simply by playing such crazy rook threats one after the other, always answered by white. Thus, white must sometimes ignore such rook moves by black. This possibility is analyzed in detail in~\cite[fig.~11]{EvansHamkins2014:TransfiniteGameValuesInInfiniteChess}, and the same analysis applies here. To summarize the conclusion there, black has no interest in moving his rooks in the rook towers except in the pattern of the main line. Black certainly gains no advantage in moving a rook in a tower to the left of the active tower, as white will simply capture in a way that activates this earlier tower, allowing white to decrease $r$, which brings about black's demise more quickly. If black attempts to amass a rook army in the outer towers, beyond the active tower---these towers are not used at all in the main line---then white can safely ignore a certain fixed proportion of black's moves, always responding when the same rook moves twice and never ignoring rook movements in adjacent towers. In this way, white will decrease $h$ more quickly than he otherwise would if black had instead played threats in the bishop gun battery, since those moves lead to much longer sequences of forced play for black to impose on white. So black is better off with the main line than attempting these irrelevant rook threats.

Next, let us consider the possibility that black attempts to find a quicker checkmate after escaping from the gateway, with the white guard two steps behind. For example, rather than going straight to the throne room, can black release a rook from the rook tower in such a way that the throne room becomes corrupted and no longer able to hold the black king? Black could use a black-square black bishop to capture at {\tt 1\ldots Bxc17} in figure~\ref{Figure:RookTowers}. But there seems to be no follow-up to this, since the rook will be captured if it moves, and in any case, white can proceed to a quick checkmate with the white guard. Since the checkmate comes so soon upon the gateway exits of the bishops, there isn't time for any complicated corruption plan.

Suppose next that black attempts to avoid loss by capturing at the outset the key mating bishop (blood-red in figure~\ref{Figure:RookTowers}), with {\tt 1\ldots Rxd15}. White simply recaptures with {\tt 2.exd15}, releasing the white bishop at {\tt f15} into the open area, with a quick checkmate now available.

Lastly, let us consider the possibility that black tries to make use of his free pawn in the pawn diagonal of the rook tower that is freed up when white captures the {\tt f7} bishop after the first move. For example, following the main line, the game may begin with {\tt 1\ldots Bq18 2.rxq18}, and then black may consider the pawn at {\tt r18}. Can he use it profitably? No. If he descends, then white will capture and poke a hole in the lower pawn diagonal, with {\tt 2\ldots r17 3.sxr17}. White can now advance the white pawn at {\tt q18} as in the main line to release a white-square white bishop into the channel above the pawn diagonal, placing it at {\tt q18}, and then advancing the white pawn on {\tt r17}, enabling a white bishop to escape into the free area and deliver checkmate in the throne room. So after white captures the black bishop on the first move, black should not move the corresponding black pawn at {\tt r18}. Incidentally, this analysis is precisely why in our current position, we have a double pawn wall as the lower diagonal for the rook towers, whereas in the corresponding position of~\cite{EvansHamkins2014:TransfiniteGameValuesInInfiniteChess}, there is only a single wall; in that position, it didn't matter if white bishops were released outside the tower position, since the checkmate configuration was different; but here, it does matter, and the double pawn wall prevents the release.

Let us now consider alternative lines of play for white. We claim that white has no quicker line to checkmate than what we have described in the main line. Unless the position were to become corrupted in some way, which we shall argue does not happen, then white must aim to release a white bishop into the open area, in order to deliver checkmate in the throne room. The only white bishops on the lower portion of the board are the white guard bishops in the gateway terminals, but these are confined to four squares each, with no possibility of escape, unless black should attack, capture and exit the doorway, in which case white will follow two steps behind as in the main line. There are no threats that white can make with those bishops, since as we have argued, black will not actually ever exit the doorway. The other white bishops are all locked up in the rook towers, and no white bishop in a tower to the left of the active tower can move.

Next, consider the possibilities for white arising from the free white-square white bishops that will accumulate in the rook towers during the normal course of the main line. Specifically, after the opening moves {\tt 1\ldots Bq18 2.rxq18, 3.q19, 4.qxp20, 5.pxq21, 6.p20,} {\tt 7.Bq18,\hfil 8.Bpo20,\hfil 9.B22p21,\hfil 10.o22,\hfil 11.o23,\hfil 12.oxn24 Rm{\it j}\break 13.lxm{\it j}}, which we discussed in connection with figure~\ref{Figure:RookTowers}, several white bishops are free to descend into the long diagonal region above the pawn diagonal at the bottom of the rook tower position. Can these bishops serve any purpose for white? No, we claim that they are useless for white, since they are isolated in this region above the pawn diagonal, and cannot leak out into the free region (unless black should cooperate by moving his black pawn down, as we discussed two paragraphs above; but we argued that black will not do that).

Consider next the possibility that white chooses to open up additional rook towers to the right of the active rook tower. That is, after the line mentioned in the previous paragraph, white is actually free to make additional captures with the white pawn via {\tt 14.qxr22, 15.rxs23} and so on. We claim that this alternative line serves no purpose for white. At best, what happens is that additional white bishops are released from the later towers into the long diagonal region just above the pawn diagonal, but as we have already noted, for white to amass an army of white bishops there offers no advantage.

We have already noted that white's responses to black's main-line threats in the bishop battery configuration are forced.

Next, let us consider whether white can profitably deviate from the main line by moving the white bishop guards or the white guard pawns into their guarding positions prophylactically, before black makes a threat. This is not helpful to white, as black has infinitely many cannons and infinitely many gateways to choose from, and so when the moment comes to fire a new cannon, he can choose a cannon where no such movements have yet taken place, and with every new bishop firing out of a cannon, he can place it so as to attack a gateway that has not yet been prematurely protected.

Consider now the possibility that white might gain an advantage by using the white guard pawns to break through the cannon wall. That is, in figure~\ref{Figure:BishopCannon}, consider that white advances one of the guard pawns prematurely, not just to guard the corresponding pawn in the cannon wall, but to attack the black pawn above it. For instance, consider {\tt 1.kxj8} after advancing the {\tt k}-pawn in figure
\ref{Figure:BishopCannon}. This is very bad for white, since black will recapture with {\tt 1\ldots kxj8}, after which black will quickly mate by escaping with his bishop on {\tt j10} or, in case that bishop has been fired, black will mate with the backup bishop on {\tt ill}. Since any breach in the cannon wall would result in black being able to escape with black-square black bishops into the free region, where the bishops would quickly deliver checkmate in the throne room, white does not want to allow or cause the wall to be breached.

In conclusion, we claim that neither player gains advantage from deviating from the main line of play, and as a result, the position is a win for white with game value $\omega^4$.

\section{A more compact version of the position}

\begin{figure}[h]
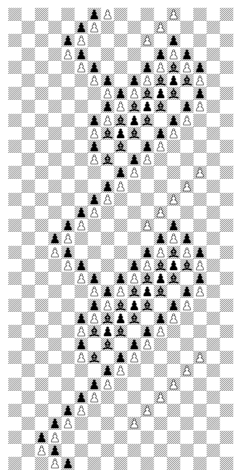

\catcode`\.=\active\def.{{1}}
\chessboard[maxfield=q35,
            boardfontsize=5pt,
            label=false,
            labelbottom=false,
            labelleft=false,
            labelfontsize=5pt,
            labelleftwidth=2ex,
            hlabelformat={\hfill\arabic{ranklabel}\,\,},
            showmover=false,
            border=false,
            margin=false,
            setfen=%
......pP....P..../%
.....pP....P....../%
....pP....P.p...../%
....Pp.....pPp..../%
.....Pp...pPbPp..../%
......Pp.pPbpbP/%
.......PpPbpb.p/%
.......pPbpb.pP../%
......pPbpb.pP../%
......Pbpb.pP.../%
......p.b.pP.../%
......Pb.pP...../%
........pP....P./%
.......pP....P/%
......pP....P/%
.....pP....P./%
....pP....P.p/%
...pP......pPp...../%
...Pp.....pPbPp.../%
....Pp...pPbpbP..../%
.....Pp.pPbpb.p/%
......PpPbpb.pP/%
......pPbpb.pP../%
.....pPbpb.pP../%
.....Pbpb.pP.../%
.....p.b.pP.../%
.....Pb.pP....P./%
.......pP....P./%
......pP....P/%
.....pP....P/%
....pP....P./%
...pP....P/%
..pP...../%
..Pp...../%
...Pp..../%
]%
\caption{Compact bishop cannon battery, with vertical right boundary}\label{Figure:CompactBishopCannon}
\end{figure}

We should like to conclude this paper by mentioning that it is possible to fit a modified version of the bishop cannon battery on a quarter-board, rather than a half board, thereby making available a large amount of space for future improvements. The idea is to retain the gateway terminals as in the main position shown in figure~\ref{Figure:MainPositionValueOmega^4}, but replace the cannons on the lower right of the position with a more compact geometry, indicated in figure~\ref{Figure:CompactBishopCannon}. The point here is that one can still have one cannon of each finite size, but in the new arrangement the right-most boundary does not change; rather, the back ends of the cannons line up, and they each extend a little into the arena, without sacrificing the possibility that each bishop fired from any of the cannons can attack arbitrarily large gateway terminals (although in this arrangement, very large cannons lose the ability to attack small gateway terminals). In this way, the fundamentals of the position will be essentially the same as in the main position of this article, with game value $\omega^4$, yet the lower right quadrant of the board will become completely empty.

We can imagine that such an arrangement might help to construct a position in infinite chess with game value $\omega^5$. Specifically, our idea is that one might augment the current position with an additional battery of cannons of some kind at the lower right, and modify the current cannons to require a capture of some kind from those cannons, with forced reply, in order to become initiated; in this way, one could fire the more remote cannons, which would allow you to fire arbitrarily large $n$ number of the current cannons successively, each time with force, and this would make an $\omega^3$ level of delay, which would make $\omega^5$ for the position overall. But we have as yet no specific actual position to propose beyond this idea.

\bibliographystyle{alpha}
\bibliography{HamkinsBiblio,MathBiblio,WebPosts}

\end{document}